UDC 004.942

**Yakiv O. Kalinovsky,** Dr.Sc., Senior Researcher,
Institute for Information Recording National Academy of Science of Ukraine, Kyiv, Shpaka str. 2, 03113, Ukraine, E-mail: kalinovsky@i.ua

**Dmitry V. Lande,** Dr.Sc., Head of Department,
Institute for Information Recording National Academy of Science of Ukraine, Kyiv, Shpaka str. 2, 03113, Ukraine, E-mail: dwlande@gmail.com

**Yuliya E. Boyarinova,** PhD, Associate Professor,
National Technical University of Ukraine "KPI", Kyiv, Peremogy av. 37, 03056, Ukraine, E-mail: ub@ua.fm

**Alina S. Turenko**, postgraduate student,
Institute for Information Recording National Academy of Science of Ukraine, Kyiv, Shpaka str. 2, 03113, Ukraine, E-mail: asturenko@mail.ru

# Clifford Type Algebra Characteristics Investigation

The main properties of hypercomplex generalization of quaternion system as antiquaternion are presented in this article. Definitions and studied of antiquaternions conjugation are introduced, their norm and zero divisor, and how to perform operations on them.

**Introduction**

There are many applications of hypercomplex number systems in the methods of processing and presenting information. The quaternion system has particular importance, applications which can solve many practical problems: the navigation and management of mobile objects, in mechanics, electrodynamics, cryptography, digital signal processing, and others.

Numerous applications of quaternions is caused by their properties which allow to carry out various operations with vectors in three-dimensional the Cartesian system of coordinates. Therefore, it is advisable to consider the properties of other hypercomplex number systems, such as a hypercomplex generalization of quaternion as antiquaternion system. Such studies will allow to solve new practical problems or facilitate decision which has been previously discussed.

**Definition of the problem**

This work investigates the question of synthesis of hypercomplex number system of antiquaternion as a result of the application to the system of complex numbers doubling procedure of Grassmann - Clifford by system of double numbers. The main properties of antiquaternions, algorithms of performance of a set of the algebraic operations, necessary for use the system of antiquaternions in mathematical modeling, are studied.

**Definitions and basic properties of antiquaternions**

It is known [1, 2], the four-dimensional hypercomplex number system called a system of quaternions $H$ with basis $\{e_1, e_2, e_3, e_4\}$ the multiplication table of elements which is;

| $H$   | $e_1$ | $e_2$  | $e_3$  | $e_4$  |
|-------|-------|--------|--------|--------|
| $e_1$ | $e_1$ | $e_2$  | $e_3$  | $e_4$  |
| $e_2$ | $e_2$ | $-e_1$ | $e_4$  | $-e_3$ |
| $e_3$ | $e_3$ | $-e_4$ | $-e_1$ | $e_2$  |
| $e_4$ | $e_4$ | $e_3$  | $-e_2$ | $-e_1$ |

Quaternions are the result of anticommutative doubling of complex numbers $C$ by the same number system. Or, using a system of markings, introduced in [3], we can write:

$$H = D(C, C). \quad (1)$$

If we redouble the system of complex numbers $C$ by system of binary numbers $W(e,2)$ with multiplication table

| $W$   | $e_1$ | $e_2$ |
|-------|-------|-------|
| $e_1$ | $e_1$ | $e_2$ |
| $e_2$ | $e_2$ | $e_1$ |

,

we obtain a system of antiquaternions $AH$, or through operator of doubling:

$$AH = D(C(e,2), W(f,2)). \quad (2)$$

Indeed, if we take the composition of bases $\{e_1, e_2\}$ and $\{f_1, f_2\}$, we obtain a basis

$$\{e_1 f_1, e_2 f_1, e_1 f_2, e_2 f_2\}.$$

Multiplication table of the obtained hypercomplex number system is constructed using multiplying of elements of this basis. However, we believe that the basic elements of the same name are multiplied by the rules of systems $C$ and $W$. When multiplying them together remains commutative only, when at least one multiplier is $e_1$ or $f_1$. Basic elements $e_2$ and $f_2$ are multiplied anticommutatative.

$$e_2 f_2 = -f_2 e_2.$$

We will give some examples of multiplication of basic elements, taking into account these rules:

$$e_1 f_1 \cdot e_1 f_1 = e_1 e_1 \cdot f_1 f_1 = e_1 f_1$$

$$e_2 f_1 \cdot e_2 f_1 = e_2 e_2 \cdot f_1 f_1 = -e_1 f_1$$

$$e_2 f_2 \cdot e_2 f_2 = -e_2 e_2 \cdot f_2 f_2 = e_1 f_1.$$

If to rename two-symbolical names of basic elements in one-symbolical:

$$e_1 f_1 \to e_1, \ e_2 f_1 \to e_2, \ e_1 f_2 \to e_3, \ e_2 f_2 \to e_4,$$

that we will receive the basic elements multiplication table of antiquaternions system:

| AH | $e_1$ | $e_2$ | $e_3$ | $e_4$ |
|---|---|---|---|---|
| $e_1$ | $e_1$ | $e_2$ | $e_3$ | $e_4$ |
| $e_2$ | $e_2$ | $-e_1$ | $e_4$ | $-e_3$ |
| $e_3$ | $e_3$ | $-e_4$ | $e_1$ | $-e_2$ |
| $e_4$ | $e_4$ | $e_3$ | $e_2$ | $e_1$ |

(3)

The principle of basic elements multiplication is represented in fig. 1, on which basic elements $e_2, e_3, e_4$ of antiquaternion system are represented by triangle tops. Product of any two elements from this three is equal to the third if movement from the first to the second multiplier coincides with the arrow direction, if the movement is opposite to the direction of the arrow - the third with a minus sign.

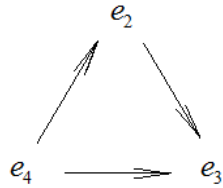

Fig. 1 Schematic image of the basic elements multiplication table of antiquaternion system.

Thus, antiquaternions are numbers of a look

$$w = a_1 e_1 + a_2 e_2 + a_3 e_3 + a_4 e_4, \qquad (4)$$

where: $a_i \in R$.

**Addition and multiplication of antiquaternions**

In antiquaternions system it is entered addition and multiplication operations as follows:
The antiquaternion $w_3$: $w_3 = w_1 + w_2 = (a_1 + b_1)e_1 + (a_2 + b_2)e_2 + (a_3 + b_3)e_3 + (a_4 + b_4)e_4$
is called *the sum* of two antiquaternions $w_1 = a_1 e_1 + a_2 e_2 + a_3 e_3 + a_4 e_4$ and $w_2 = b_1 e_1 + b_2 e_2 + b_3 e_3 + b_4 e_4$.

The antiquaternion $w_3$:

$$\begin{aligned} w_3 = w_1 w_2 &= (a_1 b_1 - a_2 b_2 + a_3 b_3 + a_4 b_4)e_1 + (a_1 b_2 + a_2 b_1 - a_3 b_4 + a_4 b_3)e_2 + \\ &+ (a_1 b_3 + a_3 b_1 - a_2 b_4 + a_4 b_2)e_3 + (a_1 b_4 + a_4 b_1 + a_2 b_3 - a_3 b_2)e_4 \end{aligned} \qquad (5)$$

is called *the product* of two antiquaternions $w_1 = a_1 e_1 + a_2 e_2 + a_3 e_3 + a_4 e_4$ and $w_2 = b_1 e_1 + b_2 e_2 + b_3 e_3 + b_4 e_4$.

According to rules of addition and multiplication of antiquaternions it is possible to mark out their main properties:
1) addition operation is commutative: $w_1 + w_2 = w_2 + w_1$;
2) addition operation is associative: $(w_1 + w_2) + w_3 = w_1 + (w_2 + w_3)$;
3) multiplication operation is noncommutative: $w_1 w_2 \neq w_2 w_1$. (6)

Really:

$$\begin{aligned} w_1 w_2 &= (a_1 e_1 + a_2 e_2 + a_3 e_3 + a_4 e_4)(b_1 e_1 + b_2 e_2 + b_3 e_3 + b_4 e_4) = \\ &= (a_1 b_1 - a_2 b_2 + a_3 b_3 + a_4 b_4)e_1 + (a_1 b_2 + a_2 b_1 - a_3 b_4 + a_4 b_3)e_2 + , \\ &+ (a_1 b_3 - a_2 b_4 + a_3 b_1 + a_4 b_2)e_3 + (a_1 b_4 + a_2 b_3 - a_3 b_2 + a_4 b_1)e_4 \end{aligned}$$

but opposite order is such as:

$$w_2 w_1 = (b_1 e_1 + b_2 e_2 + b_3 e_3 + b_4 e_4)(a_1 e_1 + a_2 e_2 + a_3 e_3 + a_4 e_4) =$$
$$= (b_1 a_1 - b_2 a_2 + b_3 a_3 + b_4 a_4)e_1 + (b_1 a_2 + b_2 a_1 - b_3 a_4 + b_4 a_3)e_2 +$$
$$+ (b_1 a_3 - b_2 a_4 + b_3 a_1 + b_4 a_2)e_3 + (b_1 a_4 + b_2 a_3 - b_3 a_2 + b_4 a_1)e_4 =$$
$$= (a_1 b_1 - a_2 b_2 + a_3 b_3 + a_4 b_4)e_1 + (a_1 b_2 + a_2 b_1 + a_3 b_4 - a_4 b_3)e_2 +$$
$$+ (a_1 b_3 + a_2 b_4 + a_3 b_1 - a_4 b_2)e_3 + (a_1 b_4 - a_2 b_3 + a_3 b_2 + a_4 b_1)e_4 \neq w_1 w_2$$

That is carried out (6).

4) multiplication operation is associative: $w_1(w_2 w_3) = (w_1 w_2)w_3$.

It can be proved directly, using (5).

5) In the same way it is possible to prove the distributivity of antiquaternions:

$$w_1(w_2 + w_3) = w_1 w_2 + w_1 w_3;$$

6) for antiquaternions is determined action of multiplication by a scalar: $k \in R$,

$$kw_1 = ka_1 e_1 + ka_2 e_2 + ka_3 e_3 + ka_4 e_4;$$

7) for $\forall k_1, k_2 \in R$ is performed $(k_1 w_1)(k_2 w_2) = k_1 k_2 (w_1 w_2)$.

**Antiquaternions norm definition**

In the work [2] the norm of hypercomplex number generally is determined by a formula

$$N(w) = \sum_{i=1}^{n} \gamma_{ij}^k a_i, \tag{7}$$

where $\gamma_{ij}^k$ - structural constants of hypercomplex number system of antiquaternions $AH$, which are defined from (3). On this basis the norm of matrix is constructed [2].

$$N(w) = \begin{vmatrix} a_1 & -a_2 & a_3 & a_4 \\ a_2 & a_1 & a_4 & -a_3 \\ a_3 & a_4 & a_1 & -a_2 \\ a_4 & -a_3 & a_2 & a_1 \end{vmatrix}. \tag{8}$$

Having calculated the matrix (8) determinant we will receive norm of hypercomplex number $w$:

$$N(w) = \left(a_1^2 + a_2^2 - a_3^2 - a_4^2\right)^2 \tag{9}$$

By analogy to the theory of quaternions we will call a root of the norm a pseudonorm of antiquaternions (9), which will be denoted as $N(w)$:

$$N(w) = a_1^2 + a_2^2 - a_3^2 - a_4^2. \tag{10}$$

Apparently from (10), the pseudonorm can be negative. It is possible to show that the pseudonorm entered by such method is multiplicative:

$$N(w_1 w_2) = N(w_1) N(w_2). \tag{11}$$

**Definition and characteristics of conjugate antiquaternions**

We introduce the antiquaternion conjugate definition of

$$\overline{w} = b_1 e_1 + b_2 e_2 + b_3 e_3 + b_4 e_4 \tag{12}$$

on the basis of equality

$$w\overline{w} = N(w), \tag{13}$$

as it is offered in [2]. If (13) substitute (5) and (10), and to equate coefficients at identical basic elements that we will receive linear algebraic system concerning variables $b_1, b_2, b_3, b_4$:

$$\begin{cases} a_1 b_1 - a_2 b_2 + a_3 b_3 + a_4 b_4 = a_1^2 + a_2^2 - a_3^2 - a_4^2 \\ a_1 b_2 + a_2 b_1 - a_3 b_4 + a_4 b_3 = 0 \\ a_1 b_3 + a_3 b_1 - a_2 b_4 + a_4 b_2 = 0 \\ a_1 b_4 + a_4 b_1 + a_2 b_3 - a_3 b_2 = 0 \end{cases}, \tag{14}$$

which solutions have the form:

$$b_1 = a_1, b_2 = -a_2, b_3 = -a_3, b_4 = -a_4 \tag{15}$$

Therefore, if the original antiquaternion $w = a_1 e_1 + a_2 e_2 + a_3 e_3 + a_4 e_4$, that the conjugate antiquaternion to it has view:

$$\overline{w} = a_1 e_1 - a_2 e_2 - a_3 e_3 - a_4 e_4. \tag{16}$$

We will define some properties of conjugate antiquaternions.
1) the sum and the product of conjugate antiquaternions is a real numbers;
2) the conjugate of the sum is the sum of conjugated $\overline{w_1 + w_2} = \overline{w_1} + \overline{w_2}$;
3) the conjugate of the product is the product of conjugated $\overline{w_1 w_2} = \overline{w_1}\,\overline{w_2}$, which can be verified directly.

**Zero divisors and their properties**

Not equal to zero antiquaternion $w_1 \neq 0$ is called as zero divider if there is such antiquaternion $w_2 \neq 0$, that their product is equal to zero $w_1 w_2 = 0$, and it means the same ratio between their pseudonorm:

$$N(w_1 w_2) = 0. \tag{17}$$

On the basis of (11) pseudonorm of a divider of zero is equal to zero

$$N(w_1) = 0. \tag{18}$$

But from (17) it follows $w_2 = \overline{w_1}$, that is, if $w \in AH$ - zero divisor, then $\overline{w}$ - zero divisor too.

From (18) follows the sign of zero divisor in this hypercomplex number system $AH$:

$$a_1^2 + a_2^2 = a_3^2 + a_4^2. \tag{19}$$

**Antiquaternions division operation**

The share from the left division of an antiquaternion $w_1$ into an antiquaternion $w_2$ is the solution of the equation

$$w_2 x = w_1 \tag{20}$$

To solve the equation (20) it is necessary to multiply its both parts at the left at first on $\overline{w_2}$, and then on $\dfrac{1}{|w_2|^2}$, where $|w_2|^2 \neq 0$. We will receive

$$x_l = \dfrac{1}{|w_2|^2} \overline{w_2} w_1. \tag{21}$$

Direct substitution of (21) in the equation (20) we find out that this expression is the solution of this equation.

We introduce the right division on the basis the equation

$$xw_2 = w_1, \tag{22}$$

from where

$$x_r = \dfrac{1}{|w_2|^2} w_1 \overline{w_2}. \tag{23}$$

Since product of antiquaternions depends on the order of the factors, then $x_l \neq x_r$. Thus, the solution of the equation (20) is called the left share, and the equation (22) - the right share [3].

It should be noted that the operation of division, unlike the fields of real and complex numbers, not possible, not only by zero but also by zero divisors in this system $AH$.

**Geometrical sense of antiquaternions**

Vector parts of antiquaternions form three-dimensional linear space which we will call imaginary space of antiquaternions. We will represent it in three-dimensional Euclidean space.

We will consider $\forall w \in AH$. Let $w = a_1 e_1 + a_2 e_2 + a_3 e_3 + a_4 e_4$ and $p$ — its pseudonorm $p = N(w)$.

We will fix scalar part.
Then

$$a_2^2 - a_3^2 - a_4^2 = p - a_1^2 \Rightarrow a_3^2 + a_4^2 - a_2^2 = a_1^2 - p.$$

We will consider geometrical sense of the given expression in three-dimensional imagined space of antiquaternions (fig. 2):

1) if $a_1^2 - p = 0$, then the set of antiquaternions forms a cone

$$a_3^2 + a_4^2 - a_2^2 = 0;$$

2) if $a_1^2 - p > 0$, then the set of antiquaternions forms a hyperboloid of one sheet

$$\dfrac{a_3^2}{a_1^2 - p} + \dfrac{a_4^2}{a_1^2 - p} - \dfrac{a_2^2}{a_1^2 - p} = 1;$$

3) if $a_1^2 - p < 0$, then set of antiquaternions in the form $w = a_1 e_1 + a_2 e_2 + a_3 e_3 + a_4 e_4$, where $a_1^2 < p$ forms a two sheeted hyperboloid

$$\dfrac{a_3^2}{a_1^2 - p} + \dfrac{a_4^2}{a_1^2 - p} - \dfrac{a_2^2}{a_1^2 - p} = -1.$$

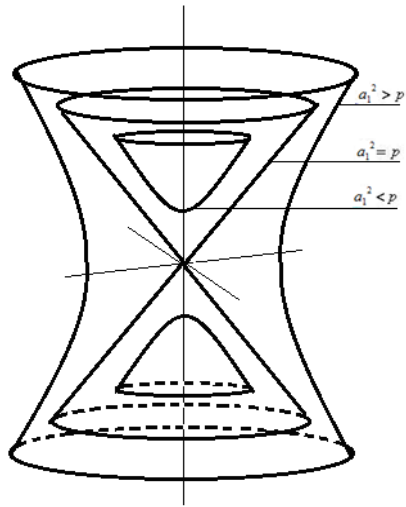

Fig. 2 Geometrical sense of antiquaternions

**Conclusions**

From the aforesaid follows that a set of arithmetic and algebraic operations in hypercomplex number system of antiquaternions is defined, which allows to use this number system for mathematical models creation in various areas of science and technology.

These operations allow to build various functions from antiquaternions, such as exponential, logarithmic, trigonometrical and hyperbolic functions, that will be a subject of further scientific researches.